\documentclass[psamsfonts]{amsart}

\usepackage{amssymb,amsfonts}
\usepackage[all,arc]{xy}
\usepackage{enumerate}
\usepackage{mathrsfs}
\usepackage{cite}

\newtheorem{theorem}{Theorem}[section]
\newtheorem{corollary}[theorem]{Corollary}
\newtheorem{proposition}[theorem]{Proposition}

\theoremstyle{definition}
\newtheorem{definition}[theorem]{Definition}

\theoremstyle{remark}
\newtheorem{remark}[theorem]{Remark}

\let\phi=\varphi

\def\Ann{\operatorname{Ann}}

\def\Spec{\operatorname{Spec}}

\let\oldbigwedge\bigwedge
\def\BIGwedge{{\textstyle\oldbigwedge}}
\def\medwedge{{\scriptstyle\oldbigwedge}}
\def\bigwedge{\mathchoice{\BIGwedge}{\BIGwedge}{\medwedge}{}}

\DeclareMathOperator{\Id}{Id}
\DeclareMathOperator{\Max}{Max}

\DeclareMathOperator{\MC}{MC}

\makeatletter

\let\c@equation\c@theorem
\makeatother
\numberwithin{equation}{section}

\begin{document}
\title[Ideals of Semirings]{Some Remarks on Ideals of Commutative Semirings}

\author[P. Nasehpour]{\bfseries Peyman Nasehpour}

\address{ 
Department of Engineering Science \\ 
Golpayegan University of Technology   \\ 
Golpayegan\\
Iran}
\email{nasehpour@gut.ac.ir, nasehpour@gmail.com}


\begin{abstract}

The main purpose of this paper is to investigate prime, primary, and maximal ideals of semirings. The localization and primary decomposition of ideals in semirings are also studied. 
\end{abstract}

\maketitle

\section{Introduction}

Semirings are ring-like algebraic structures that subtraction is either impossible or disallowed. Other ring-like algebraic structures include pre-semirings \cite{GondranMinoux2008}, hemirings \cite{Golan1999(b)}, and near-rings \cite{Pilz1983}.

Vandiver introduced a simple type of a ring-like algebraic structure, with the term ``semi-ring'', in which the cancellation law of addition does not hold \cite{Vandiver1934}. In many references (see the explanations in page 3 of the book \cite{Glazek2002}), a semiring is an algebraic structure $(S,+,\cdot,0,1)$ with the following properties:

\begin{enumerate}
	\item $(S,+,0)$ is a commutative monoid,
	\item $(S,\cdot,1)$ is a monoid with $1\neq 0$,
	\item $a(b+c) = ab+ac$ and $(b+c)a = ba+ca$ for all $a,b,c\in S$,
	\item $a\cdot 0 = 0\cdot a = 0$ for all $a\in S$. 
\end{enumerate}

Note that while the last axiom, i.e. $a\cdot 0 = 0\cdot a = 0$ for all $a\in S$, is omitted from the definition of a ring, since it follows from the other ring axioms, but here it does not, and it is necessary to state it in the definition (see Example 5.1.2. in \cite{GondranMinoux2008}). A semiring $S$ is commutative if $ab = ba$ for all $a,b\in S$. In this paper all semirings are commutative.

Semirings have applications in science and engineering especially in computer science and are an interesting generalization of rings and bounded distributive lattices \cite{Golan1999(b)}. They can also be used to model algebraic properties of probability and modular functions \cite{NasehpourParvardi2018}. For general books on semiring theory, one may refer to the resources \cite{AhsanMordesonShabir2012, Eilenberg1974, Golan1999(b), Golan2003, GondranMinoux2008, HebischWeinert1998, KuichSalomaa1986}. 

The ideal theoretic method for studying commutative rings has a long and fruitful history \cite{Huckaba1988}. Some of the topics related to the ideals of commutative rings have been generalized and investigated for semirings \cite{GhalandarzadehNasehpourRazavi2017, Iseki1956, LaGrassa1995, Lescot2015, Nasehpour2016, Nasehpour2017, NasehpourP2018, NasehpourPP2018, Nasehpour2018}. Also, see Chapter 7 of the book \cite{Golan1999(b)}. The main purpose of this paper is to investigate the prime, primary, and maximal ideals of semirings and related concepts such as localization and primary decomposition. Here is a brief sketch of the contents of our paper:

In Section \ref{sec:ideals}, we bring some primitive properties of ideals in semirings. We will use those properties in the paper, sometimes without mentioning them. We also fix some terminologies in this section. In Section \ref{sec:primes}, we investigate prime and maximal ideals of semirings. Let us recall that $W$ is an MC-set if and only if it is a submonoid of $(S,\cdot)$. In this section, similar to commutative algebra, we prove that the maximal elements of the set of all ideals disjoint from an MC-set of a semiring are prime ideals (check Theorem \ref{maxisprime1}). Note that in the proof of this theorem, we use ``Zorn's lemma'' and since one of the corollaries of this important theorem is that any semiring has a maximal ideal (Corollary  \ref{maxisprime2}), one may ask if the converse holds, i.e. if each semiring has at least one maximal ideal, then Zorn's lemma is true. In fact, in Corollary \ref{Hodges}, by using the main theorem of a paper by Hodges \cite{Hodges1979}, we see that in Zermelo-Fraenkel set theory, the following statements are equivalent:

\begin{enumerate}
	\item Axiom of Choice holds;
	\item Every semiring with $1\neq 0$ has a maximal ideal;
	\item Every commutative ring with $1\neq 0$ has a maximal ideal;
	\item Every unique factorization domain has a maximal ideal.
\end{enumerate} 

Let us recall that a semiring $S$ is called \emph{Artinian} if any descending chain of ideals of $S$ stabilizes, i.e., if $I_1 \supseteq I_2 \supseteq \cdots \supseteq I_n \supseteq \cdots$ is a descending chain of ideals of $S$, then there is an $m\in \mathbb N$ such that $I_m = I_{m+k}$ for all $k\geq 0$.
In this section, we also prove that if $S$ is Artinian, it has finitely many maximal ideals (see Theorem \ref{artinianisemilocal}). 

Section \ref{sec:CE} is devoted to contraction and extension of ideals in semirings. We use some of the general results of this section for proving the statements for localization of semirings and semimodules given in Section \ref{sec:fractions}. For instance in Section \ref{sec:fractions}, we prove that if $S$ is a semiring and $U$ an MC-set of $S$, then every ideal of $S_U$ - localization of $S$ at $U$ - is an extended ideal. Also, in Theorem \ref{primesoflocalization}, we show that if $S$ is a semiring and $U$ an MC-set of $S$, then the prime ideals of $S_U$ are in one-to-one correspondence with the prime ideals of $S$ disjoint from $U$.

In Section \ref{sec:primary}, we bring the primitive properties of primary ideals of semirings. Then we pass to Section \ref{sec:decompositions} to discuss irreducible ideals and primary decompositions. Let us recall that an ideal $I$ of a semiring $S$ is called irreducible if for any ideals $J$ and $K$ of $S$, $I = J \cap K$ implies that $I = J$ or $I = K$. We also recall that an ideal $J$ of a semiring $S$ is subtractive if $a+b\in J$ and $a\in J$ implies that $b\in J$ for all $a,b\in S$. In Theorem \ref{subirreprimary}, we show that if $S$ is a Noetherian semiring and $I$ a subtractive irreducible ideal of $S$, then it is primary.

Finally, let us recall that a semiring $S$ is Noetherian if and only if each ideal of $S$ is finitely generated \cite[Proposition 6.16]{Golan1999(b)}. A semiring is subtractive if each ideal of $S$ is subtractive. In Corollary \ref{primarydecomposition}, we prove that if $S$ is a subtractive Noetherian semiring, then each ideal of $S$ can be represented as an intersection of a finite number of primary ideals of $S$.

\section{Ideals of semirings and operations on them}\label{sec:ideals}

The concept of ideals for semirings was introduced by Bourne in \cite{Bourne1951}.

\begin{definition}
	
	\label{ideal-def}
	
	A nonempty subset $I$ of a semiring $S$ is said to be an \emph{ideal} of $S$, if $a+b \in I$ for all $a,b \in I$ and $sa \in I$ for all $s \in S$ and $a \in I$. It is clear that the zero element $0$ belongs to any ideal of $S$. We denote the set of all ideals of $S$ by $\Id(S)$. An ideal $I$ of a semiring $S$ is called a \emph{proper ideal} of the semiring $S$ if $I \neq S$.
	
\end{definition}

\begin{proposition}
	Let $S$ be a semiring and $I$, $J$, and $K$ be ideals of $S$. If we define the addition and multiplications as follows: \[I+J := \{x+y : x\in I, y\in J\} \text{~and~} I\cdot J := \{\sum_{i\leq k} x_i y_i : x_i \in I, y_i \in J, k\in \mathbb N \}, \] then the following statements hold:
	
	\begin{enumerate}
		
		\item The sets $I+J$ and $I\cdot J$ are ideals of $S$. 
		\item $I+(J+K) = (I+J)+K$ and $I(JK) = (IJ)K$.
		\item $I+J = J+I$ and $ IJ = JI $.
		\item $I(J+K) = IJ + IK$.
		\item $I+I = I$, $I+(0) = I$, $IS=I$, $I(0)=(0)$ and $I+S = S$.
		\item If $I+J = (0)$ then $I = J = (0)$.
		\item $IJ \subseteq I \cap J$ and if $I+J = S$, then $IJ = I \cap J$.
		\item $(I+J) (I \cap J) \subseteq IJ$
		
	\end{enumerate}
	
\end{proposition}

It is easy to verify that if $I_{\alpha}$ is a family of ideals of the semiring $S$, then $\bigcap_{\alpha} I_{\alpha}$ is also an ideal of $S$.

Note that if $\{I_{\alpha}\} _{\alpha \in A}$ is a family of ideals of $S$, it can be easily seen that the subset $\{ s_{\alpha_1}+ s_{\alpha_2}+ \cdots + s_{\alpha_n} : s_{\alpha_j}\in I_{\alpha_j}, \alpha_j \in A, n\in \mathbb N \}$ of $S$, denoted by $\sum_{\alpha \in A} I_{\alpha}$, is an ideal of $S$. More generally if $A$ is a subset of a semiring $S$, we denote the set of the intersection of all ideals of $S$, which contain $A$ by $(A)$. Since $S$ contains $A$, this intersection is nonempty. Obviously, $(A)$ is the smallest ideal containing $A$. One can easily see that the elements of $(A)$ can be obtained from all possible linear combinations of elements of $A$. This is perhaps why the ideal $(A)$ is said to be an ideal of $S$ generated by $A$. Note that if $\{I_{\alpha}\} _{\alpha \in A}$ is a family of ideals of $S$, then $\sum_{\alpha \in A} I_{\alpha}$ is generated by $\bigcup I_{\alpha}$.

Let us recall that in semiring theory, the multiplication of ideals distributes over finite addition of ideals. Also, it is a routine exercise to check that if $J$ is also an ideal of $S$, then $J \cdot \sum_{\alpha \in A} I_{\alpha} = \sum_{\alpha \in A} J \cdot I_{\alpha}$, which means that multiplication of ideals distributes over infinite addition of ideals as well. One may interpret some of these properties as follows: 

\begin{proposition}
	Let $S$ be a semiring. If we denote the set of all ideals of $S$ by $\Id(S)$, then the following statements hold:

	\begin{enumerate}
		\item $(\Id(S),+,\cdot)$ is an additively-idempotent semiring.
		\item $(\Id(S),\subseteq)$ is a bounded lattice, where $(0)$ is the least and $S$ is the greatest element of this lattice and $\sup\{I,J\} = I+J$ and $\inf\{I,J\} = I \cap J$.
		\item $(\Id(S), \Sigma ,\cdot)$ is an example of $c$-semirings \cite[2.1 p. 23]{Bistarelli2004}.
	\end{enumerate}
\end{proposition}

For any ideals $I,J$ of a semiring $S$, it is defined that $[I \colon J]=\lbrace s\in S \colon sJ \subseteq I \rbrace$.

\begin{proposition}
	
	\label{ideal:ideal}
	
	Let $I,J,K, I_{\alpha},$ and $J_{\alpha}$ be ideals of a semiring $S$. The following statements hold:
	
	\begin{enumerate}
		
		\item $I \subseteq [I : J]$
		\item $[I:J]J \subseteq I$
		\item $[[I : J] : K] = [I:JK] = [[I : K] : J]$
		\item $[\bigcap_{\alpha} I_{\alpha}: J] = \bigcap_{\alpha} [I_{\alpha}: J]$
		\item $[I : \sum_{\alpha} J_{\alpha}] = \bigcap_{\alpha} [I : J_{\alpha}]$
		\item $[I : J] = [I : I+J]$
		
	\end{enumerate}
	
\end{proposition}

Cancellation ideals were introduced and investigated by Susan LaGrassa in her Ph.D. Thesis \cite{LaGrassa1995}.

\begin{definition}
A nonzero ideal $I$ of a semiring $S$ is called a \emph{cancellation ideal}, if $IJ=IK$ implies $J=K$ for all ideals $J$ and $K$ of $S$.
\end{definition}

\begin{remark}
Let us recall that an element $s$ of a semiring $S$ is said to be multiplicatively cancelable, if $sb=sc$ implies $b=c$ for all $b,c \in S$. If $s$ is a multiplicatively cancelable element of $S$, then the principal ideal $(s)$ is a cancellation ideal and the proof is as follows:
	
	Take $I,J$ to be arbitrary ideals of $S$ such that $(s)I = (s)J$ and imagine $x\in I$, then obviously $sx \in (s)J$, which means that there is a $y\in J$ such that $sx = sy$. But $s$ is a multiplicatively cancelable element. So $x=y$, which implies that $I \subseteq J$. Similarly it is proved that $J\subseteq I$, which means that $(s)$ is cancellation.
\end{remark}

The following proposition taken from \cite{GhalandarzadehNasehpourRazavi2017} is the semiring version of a statement mentioned in \cite[Exercise. 4, p. 66]{Gilmer1972}.

\begin{proposition}
	
	\label{can2}
	
	Let $S$ be a semiring and $I$ be a nonzero ideal of $S$. Then the following statements are equivalent:
	
	\begin{enumerate}
		
		\item $I$ is a cancellation ideal of $S$,
		
		\item $[IJ \colon I]=J$ for any ideal $J$ of $S$,
		
		\item $IJ \subseteq IK$ implies $J \subseteq K$ for all ideals $J,K$ of $S$.
		
	\end{enumerate}
	
	\begin{proof}
		By considering this point that the equality $[IJ \colon I]I=IJ$ holds for all ideals $I,J$ of $S$, it is, then, easy to see that (1) implies (2). The rest of the proof is straightforward.
	\end{proof}
	
\end{proposition}

Let $I$ be an ideal of $S$. The \emph{radical} of $I$, denoted by $\sqrt{I}$, is defined to be the set of all elements of the form $s^n$, where $s\in I$ and $n\in \mathbb N$.

\begin{proposition}
	Let $S$ be a semiring and $I,J$ be ideals of $S$. Then the following statements hold:
	
	\begin{enumerate}
		
		\item $ I \subseteq \sqrt{I}$ and $\sqrt{I} = \sqrt{\sqrt{I}}$.
		
		\item $ \sqrt{IJ} = \sqrt{I \cap J} = \sqrt{I} \cap \sqrt{J}$.
		
		\item $ \sqrt{I}=S$ if and only if $I=S$.
		
		\item $ \sqrt{I+J} = \sqrt{\sqrt{I} + \sqrt{J}}$
		
	\end{enumerate}
	
	\begin{proof}
		Straightforward.
	\end{proof}

\end{proposition}

\section{Prime and Maximal Ideals of Semirings}\label{sec:primes}

\subsection*{Prime ideals} We start this section by defining prime ideals:

\begin{definition}
	
	A proper ideal $P$ of a semiring $S$ is called a \emph{prime ideal} of $S$, if $ab\in P$ implies either $a\in P$ or $b\in P$. We denote the set of all prime ideals of $S$ by $\Spec(S)$.
	
\end{definition}

It is straightforward to see that $P$ is a prime ideal of $S$ if and only if $P \neq S$ and $IJ \subseteq P$ implies either $I \subseteq P$ or $J \subseteq P$ for all ideals $I$ and $J$ of $S$ \cite[Corollary 7.6]{Golan1999(b)}. This implies the following proposition:

\begin{proposition}
	
	\label{primecontainingideals}
	
	Let $P$ be a prime ideal of a semiring $S$ and $I_1, I_2, \ldots, I_n$ be arbitrary ideals of $S$. Then the following statements are equivalent:
	
	\begin{enumerate}
		
		\item $P \supseteq I_k$ for some $1 \leq k \leq n$.
		\item $P \supseteq \bigcap_{1 \leq k \leq n} I_k$.
		\item $P \supseteq \prod_{1 \leq k \leq n} I_k$.
		
	\end{enumerate}
	
\end{proposition}

A nonempty subset $W$ of a semiring $S$ is said to be a \emph{multiplicatively closed set} (for short an MC-set) if $1\in W$ and for all $w_1,w_2 \in W$, we have $w_1 w_2 \in W$. In other words, $W$ is an MC-set if and only if it is a submonoid of $(S,\cdot)$. It is clear that an ideal $P$ of $S$ is a prime ideal of $S$ if and only if $S-P$ is an MC-set. The following theorem is semiring version of a theorem in commutative algebra due to German mathematician Wolfgang Krull (1899-1971):

\begin{theorem}
	
	\label{maxisprime1}
	
	The maximal elements of the set of all ideals disjoint from an MC-set of a semiring are prime ideals.
	
	\begin{proof}
		Let $S$ be a semiring and $W \subseteq S$ an MC-set. Let $\sum$ be the set of all ideals disjoint from $W$. If $\{I_{\alpha}\}$ is a chain of ideals belonging to $\sum$, then $\bigcup I_{\alpha}$ is also an ideal disjoint from $W$ and an upper bound for the chain $\{I_{\alpha}\}$. Therefore according to Zorn's Lemma, $\sum$ has a maximal element. Let $P$ be a maximal element of $\sum$. We prove that $P$ is actually a prime ideal of $S$.
		
		Let $a\notin P$ and $b\notin P$. Then obviously $P\subset P+(a)$ and $P\subset P+(b)$. This means that $P+(a)$ and $P+(b)$ are ideals of $S$ such that they are not disjoint from $W$. So there exist $w_1,w_2 \in W$ such that $w_1 = p_1 + xa$ and $w_2 = p_2 + yb$ for some $p_1, p_2 \in P$ and $x,y \in S$. But $w_1 w_2 = p_1 p_2 + p_1 y b + p_2 x a + xy ab$. Now it is obvious that if $ab\in P$, then $w_1 w_2 \in P$, which contradicts this fact that $P$ is disjoint from $W$. Therefore $ab\notin P$ and $P$ is a prime ideal of $S$.
	\end{proof}
	
\end{theorem}

\subsection*{Maximal ideals} We continue this section by investigating maximal ideals.

\begin{definition}
	
	Let $S$ be a semiring. An ideal $\mathfrak{m}$ of the semiring $S$ is called a maximal ideal of $S$, if $\mathfrak{m} \subseteq I \subseteq S$ for any ideal $I$ of $S$ implies either $I = \mathfrak{m}$ or $I = S$. We denote the set of all maximal ideals of $S$ by $\Max(S)$.
	
\end{definition}

\begin{corollary}
	
	\label{maxisprime2}
	
	Any semiring $S$ possesses at least one maximal ideal and all maximal ideals of $S$ are prime ideals.
	
	\begin{proof}
		In Theorem \ref{maxisprime1}, take $W=\{1\}$.
	\end{proof}
\end{corollary}

\begin{theorem}
	
	\label{propersubsetmaximal}
	
	Any proper ideal of $S$ is a subset of a maximal ideal of $S$.
	
	\begin{proof}
		If $I$ is a proper ideal of $S$ (i.e., $I\neq S$), then a chain of all proper ideals of $S$ containing $I$ has an upper bound (the union of all those ideals) and by Zorn's Lemma, the proper ideals containing $I$ has at least one maximal element that is, in fact, a maximal ideal of $S$. This means that any proper ideal $I$ of $S$ is a subset of a maximal ideal of $S$.
	\end{proof}
	
\end{theorem}

\begin{remark}
	The proof of Theorem 1.3 in \cite{AtiyahMacDonald1969} shows that ``Axiom of Choice'' (which is equivalent to Zorn's lemma \cite[\S 3]{Monk1969}) implies that every commutative ring with $1 \neq 0$ has a maximal ideal. Dana Scott \cite{Scott1954} had asked whether the converse holds: If every commutative ring with $1 \neq 0$ has a maximal ideal, then the Axiom of Choice is true. The answer is ``yes''. In fact, Wilfrid Hodges in \cite{Hodges1979} proved the following:
	
	\begin{theorem}
		In Zermelo-Fraenkel set theory, the statement ``Every unique factorization domain has a maximal ideal'' implies the Axiom of Choice. 
	\end{theorem}
\end{remark}

\begin{corollary}
	
	\label{Hodges}
	In Zermelo-Fraenkel set theory, the following statements are equivalent:
	
	\begin{enumerate}
		\item Axiom of Choice holds;
		\item Every semiring with $1\neq 0$ has a maximal ideal;
		\item Every commutative ring with $1\neq 0$ has a maximal ideal;
		\item Every unique factorization domain has a maximal ideal.
	\end{enumerate}
\end{corollary}

\begin{remark}[Zariski Topology for Semirings]
	
	\label{zariskitopology}
	
	By Corollary \ref{maxisprime2}, it is clear that $\emptyset \neq \Max(S) \subseteq \Spec(S) \subseteq \Id(S)$. By Theorem \ref{propersubsetmaximal}, $V(I) = \{P\in \Spec(S): P \supseteq I\}$ is a nonempty set if and only if $I$ is a proper ideal of $S$. On the other hand, one can easily check that $V(I_1) \cup V(I_2) = V(I_1 \cap I_2)$ and $\bigcap_{\alpha} V(I_{\alpha}) = V(\sum_{\alpha} I_{\alpha})$. Also, $V(0) = \Spec(S)$ and $V(S) = \emptyset$. From this it follows that $\mathcal{C} = \{V(I) : I\in \Id(S)\}$ defines a topology on $\Spec(S)$, known as \emph{Zariski's topology}, which its closed sets are all elements of $\mathcal{C}$. Zariski topology, due to Kiev-born American mathematician Oscar Zariski (1899--1986), is an important topology used in algebraic geometry. Zariski with his French student, Pierre Samuel (1921--2009), wrote a two-volume book in commutative algebra \cite{ZariskiSamuel1975,ZariskiSamuel1976} that is a classic.
\end{remark}

The ring version of the following theorem is credited to Wolfgang Krull (1899-1971):

\begin{theorem}
	
	\label{krullintersectionprimes}
	
	Let $S$ be a semiring and $I$ an ideal of $S$. Then the following statements hold:
	
	\begin{enumerate}
		
		\item $\sqrt{I} = \bigcap_{P\in V(I)} P$, where $V(I) = \{P\in \Spec(S) : P \supseteq I\}$.
		
		\item $\sqrt{I}$ is an ideal of $S$.
		
	\end{enumerate}
	
	\begin{proof}
		$(1)$: It is straightforward that $\sqrt{I} \subseteq \bigcap_{P\in V(I)} P$. Now let $s\notin \sqrt{I}$. It is clear that $W_s = \{s^n: n\geq 0\}$ is an MC-set of $S$ disjoint from $\sqrt{I}$. So there exists a prime ideal containing $I$ and not containing $s$.
		
		$(2)$: Since $\sqrt{I}$ is an intersection of some ideals, it is an ideal and this completes the proof.
	\end{proof}
	
\end{theorem}

An element $s$ of a semiring $S$ is said to be \emph{invertible (unit)} if there is an $s^{\prime}\in S$ such that $s \cdot s^{\prime}=1$. The set of all invertible elements of $S$ is denoted by $U(S)$. It is obvious that $U(S)$ is an Abelian multiplicative group and is called the group of units of $S$. Obviously, $I$ is a proper ideal of $S$ if and only if it contains no invertible element of $S$.

\begin{proposition}
	
	\label{units3}
	
	Let $S$ be a semiring. Then $U(S) = S - (\bigcup_{\mathfrak{m}\in \Max(S)} \mathfrak{m})$, where by $U(S)$ we mean the set of all invertible elements of $S$.
	
	\begin{proof}
		Let $S$ be a semiring and take $U(S)$ to be the set of all invertible elements of $S$. If $s\in U(S)$, then $s$ cannot be an element of a maximal ideal of $S$. On the other hand if $s$ is not invertible, then the principal ideal $(s)$ of $S$ is proper and by Theorem \ref{propersubsetmaximal}, $(s)$ is contained in a maximal ideal $\mathfrak{m}$ of $S$ and therefore $s\in \mathfrak{m}$.
	\end{proof}
	
\end{proposition}

\begin{corollary}
	
	\label{localsemiring1}
	
	Let $S$ be a semiring. Then the following statements hold:
	
	\begin{enumerate}
		
		\item The semiring $S$ has a unique maximal ideal if and only if $S-U(S)$ is an ideal of $S$.
		
		\item The semiring $S$ is a semifield if and only if $(0)$ is a maximal ideal of $S$.
		
	\end{enumerate}
	
\end{corollary}

\begin{definition}
	
	\begin{enumerate}
		\item $(S, \mathfrak{m})$ is a local semiring if $S$ is a semiring and $\mathfrak{m}$ is its unique maximal ideal.
		
		\item A semiring $S$ is semi-local if it possesses a finite number of maximal ideals, i.e., $\mid \Max(S) \mid < \infty$.
		
		\item Two ideals $I,J$ of $S$ are called comaximal if $I+J = S$.
		
		\item The ideals $\{I_k\}^n_{k=1}$ of $S$ are said to be pairwise comaximal if $I_k + I_l = S$ for any $1 \leq k<l\leq n$.
	\end{enumerate}
	
\end{definition}

\begin{proposition}
	
	\label{comaximal}
	
	Let $S$ be a semiring. Then the following statements are equivalent.
	
	\begin{enumerate}
		\item If $I,J$ are comaximal, then $I \cap J = IJ$.
		\item If the ideals $\{I_k\}^n_{k=1}$ are pairwise comaximal, then $\bigcap^n_{k=1} I_k = \prod^n_{k=1} I_k$.
		\item If $\{ \mathfrak{m}_k \}^n_{k=1}$ is a set of $n$ distinct maximal ideals of $S$, then they are pairwise comaximal and $\bigcap^n_{k=1} \mathfrak{m}_k = \prod^n_{k=1} \mathfrak{m}_k$.
		\item $I$ and $J$ are comaximal if and only if $\sqrt{I}$ and $\sqrt{J}$ are comaximal.
	\end{enumerate}

    \begin{proof}
	
	(1): Let $I$ and $J$ be comaximal. Then $I \cap J = (I+J) (I \cap J) \subseteq IJ \subseteq I \cap J$.\\
	
	(2): Let the ideals $\{I_k\}^n_{k=1}$ be pairwise comaximal. Let $J = \bigcap^{n-1}_{k=1} I_k$. We claim that $J$ and $I_n$ are comaximal. Suppose not. Then $J+I_n$ is a proper ideal of $S$ and so is contained in a maximal ideal $\mathfrak{m}$ of $S$. This implies that $J \subseteq \mathfrak{m}$, which causes $I_k \subseteq \mathfrak{m}$ for some $1 \leq k \leq (n-1)$. From this we get that $I_k + I_n \subseteq \mathfrak{m}$, which contradicts our assumption that $I_k$ and $I_n$ are comaximal. So $J$ and $I_n$ are comaximal and $J \cap I_n = JI_n$. Now by induction, the claim $\bigcap^n_{k=1} I_k = \prod^n_{k=1} I_k$ is proved.\\
	
	(3): Let $\mathfrak{m}_1$ and $\mathfrak{m}_2$ be two distinct maximal ideals of $S$. Then clearly $\mathfrak{m}_1 \subset \mathfrak{m}_1 + \mathfrak{m}_2$ and by maximality of $\mathfrak{m}_1$, we have $\mathfrak{m}_1 + \mathfrak{m}_2 = S$. Therefore if $\mathfrak{m}_1$ and $\mathfrak{m}_2$ of $S$ are two distinct maximal ideals of $S$, $ \mathfrak{m}_1 \cap \mathfrak{m}_2 = \mathfrak{m}_1 \mathfrak{m}_2$. From this we get that if $\{ \mathfrak{m}_k \}^n_{k=1}$ is a set of $n$ distinct maximal ideals of $S$, they are pairwise comaximal and $\bigcap^n_{k=1} \mathfrak{m}_k = \prod^n_{k=1} \mathfrak{m}_k$.
	
	(4): Straightforward.
	\end{proof}
\end{proposition}

\begin{remark}
	
	\label{comaximal2}
	
	In this remark, without using Zorn's lemma, we give an alternative proof of this fact that if the ideals $\{I_k\}^n_{k=1}$ are pairwise comaximal, then $\bigcap^n_{k=1} I_k = \prod^n_{k=1} I_k$.
	
	\begin{proof}
		
		In Proposition \ref{comaximal}, we have seen that if $I+J = S$, i.e., $I,J$ are comaximal ideals of $S$, then $I \cap J = IJ$. Now fix a natural number $n>2$ and suppose that any family of pairwise comaximal ideals of $\{I_k\}^n_{k=1}$ of $S$ has this property that $\bigcap^n_{k=1} I_k = \prod^n_{k=1} I_k$. Now we prove that the statement is also true for $n+1$ pairwise comaximal ideals of $S$. Set $A = \bigcap^{n+1}_{k=1} I_k$. Now we have $A = A\cdot S = A \cdot (I_1 + \cdots + I_{n+1}) = A \cdot I_1 + \cdots + A \cdot I_{n+1}$. But by induction hypothesis, $A \cdot I_k \subseteq \prod^{n+1}_{k=1} I_k$ and the proof is complete.
	\end{proof}
	
\end{remark}

\begin{definition}
	
	A semiring $S$ is called \emph{Artinian} if any descending chain of ideals of $S$ stabilizes, i.e., if $I_1 \supseteq I_2 \supseteq \cdots \supseteq I_n \supseteq \cdots$ is a descending chain of ideals of $S$, then there is an $m\in \mathbb N$ such that $I_m = I_{m+k}$ for all $k\geq 0$.
	
\end{definition}

\begin{theorem}
	
	\label{artinianisemilocal}
	
	If $S$ is an Artinian semiring, then $S$ is semi-local.
	
	\begin{proof}
		Let $S$ be an Artinian semiring and $\{ \mathfrak{m}_k \}^{\infty}_{k=1}$ a family of infinite distinct maximal ideals of $S$. We claim that $ \mathfrak{m}_1 \mathfrak{m}_2 \cdots \mathfrak{m}_l \supset \mathfrak{m}_1 \mathfrak{m}_2 \cdots \mathfrak{m}_{l+1}$. On the contrary, if $ \mathfrak{m}_1 \mathfrak{m}_2 \cdots \mathfrak{m}_l = \mathfrak{m}_1 \mathfrak{m}_2 \cdots \mathfrak{m}_{l+1}$, then by Remark \ref{comaximal}, $J \cap \mathfrak{m}_{l+1} = J$, where $J = \mathfrak{m}_1 \mathfrak{m}_2 \cdots \mathfrak{m}_l$. This implies that $J \subseteq \mathfrak{m}_{l+1}$, which causes $\mathfrak{m}_i = \mathfrak{m}_{l+1}$ for some $1 \leq i \leq l$ that is in contradiction with the distinctness of the maximal ideals $\{ \mathfrak{m}_k \}^{\infty}_{k=1}$. This gives us the following descending chain of ideals of $S$: $ \mathfrak{m}_1 \supset \mathfrak{m}_1 \mathfrak{m}_2 \supset \cdots \supset \mathfrak{m}_1 \mathfrak{m}_2 \cdots \mathfrak{m}_l \supset \cdots$ and the proof is complete.
	\end{proof}
	
\end{theorem}

\section{Contraction and Extension of Ideals}\label{sec:CE}

First, we define homomorphism between semirings.

\begin{definition}
	
Let $S$ and $B$ be two semirings. By a semiring homomorphism from $S$ to $B$, we mean a function $\lambda: S \longrightarrow B$ with the following properties:

\begin{enumerate}
	
	\item $\lambda (r+s) = \lambda (r) + \lambda (s)$ and $\lambda (rs) = \lambda (r) \lambda (s)$ for all $r,s \in S$;
	
	\item $\lambda (0) = 0$ and $\lambda (1) = 1$.
	
\end{enumerate}

\end{definition}

\subsection*{Contraction of Ideals}

Let $S$ and $B$ be two semirings and $\varphi : S \rightarrow B$ a semiring homomorphism. If $J$ is an ideal of $B$, then $\varphi^{-1} (J)$ is an ideal of $S$ and is called contraction of $J$ and is denoted by $J^c$ or sometimes $J \cap S$. In particular, $\varphi^{-1} (0)$ is an ideal of $S$, known as the kernel of $\varphi$ and is denoted by $\ker(\varphi)$. Anyhow the kernel of a semiring homomorphism does not obey the rules of a kernel of a ring homomorphism. For example, if $\varphi$ is injective, then $\ker(\varphi)=(0)$, while the converse of this statement is not true. To check this, let $S=\{0,s,1\}$ be a totally ordered set and consider the semiring $(S,\max,\min)$ and define a function $\gamma$ from $S$ to the Boolean semiring $\mathbb B=\{0,1\}$ by $\gamma(0) = 0$ and $\gamma(s) = \gamma(1) = 1$. It is easy to see that $\gamma$ is a semiring homomorphism with $\ker(\gamma)=(0)$, while it is not one-to-one. For more on kernels of semiring homomorphisms, refer to \cite[Chap. 10]{Golan1999(b)}.

The basic properties of contraction of ideals are collected in the following.

\begin{proposition} If $J, J_1, J_2$ are ideals of a semiring $B$ and $\varphi : S \rightarrow B$ is a semiring homomorphism, then the following statements hold.
	
	\begin{enumerate}
		
		\item $(J_1 + J_2)^c \supseteq J^c_1 + J^c_2$.
		
		\item $(J_1 \cap J_2)^c = J^c_1 \cap J^c_2$.
		
		\item $(J_1 \cdot J_2)^c \supseteq J^c_1 \cdot J^c_2$.
		
		\item $(\sqrt{J})^c = \sqrt{J^c}$.
		
		\item If $Q$ is a prime ideal of $B$, then $Q^c$ is a prime ideal of $S$.
		
	\end{enumerate}
	
	\begin{proof}
		Straightforward.
	\end{proof}
	
\end{proposition}

\subsection*{Extension of Ideals}

Let $S$ and $B$ be two semirings and $\varphi : S \rightarrow B$ a semiring homomorphism. If $I$ is an ideal of $S$, then the set $\varphi(I) \subseteq B$ does not need to be an ideal of $B$. Then the extension $I^{e}$ of $I$ is defined to be the ideal generated by $\varphi(I)$ in $B$. One can easily check that $I^{e} = \{ \Sigma_{i=1}^n a_i f_i: a_i \in I, f_i \in B, i\in \mathbb N \}$. The extension of the ideal $I$ is sometimes denoted by $IB$. The basic properties of extension of ideals are collected in the following.

\begin{proposition} If $I, I_1, I_2$ are ideals of a semiring $S$ and $\varphi : S \rightarrow B$ is semiring homomorphism, then the following statements hold.
	
	\begin{enumerate}
		
		\item $(I_1 + I_2)^e = I^e_1 + I^e_2$.
		
		\item $(I_1 \cap I_2)^e \subseteq I^e_1 \cap I^e_2$.
		
		\item $(I_1 \cdot I_2)^e = I^e_1 \cdot I^e_2$.
		
		\item $(\sqrt{I})^e \subseteq \sqrt{I^e}$.
		
	\end{enumerate}
	
	\begin{proof}
		Straightforward.
	\end{proof}
	
\end{proposition}

Note that in general if $P$ is a prime ideal of $S$, then its extension $PB$ does not need to be a prime ideal of $B$. But in content semialgebras, primes extend to primes \cite[Proposition 31]{Nasehpour2016}. We end this section with the following proposition:

\begin{proposition} If $I$ is an ideal of a semiring $S$ and $J$ is an ideal of $B$ and $\varphi : S \rightarrow B$ is a semiring homomorphism, then the following statements hold.
	
	\begin{enumerate}
		
		\item $I \subseteq I^{ec}$, $J \supseteq J^{ce}$.
		
		\item $I^{e} = I^{ece}$, $J^c = J^{cec}$.
		
	\end{enumerate}
	
	\begin{proof}
		Straightforward.
	\end{proof}
	
\end{proposition}

\section{Semirings and Semimodules of Fractions}\label{sec:fractions}

Localization is a very powerful tool in commutative algebra. While apparently not all the techniques of localization are valid in commutative semiring theory, but still some of them work efficiently. In this section, we introduce the semirings and semimodules of fractions that is nothing but the localization of these algebraic objects.

Let $S$ be a semiring and $U \subseteq S$ an MC-set. Define $\sim$ on $S \times U$ by $(x,u) \sim (y,v)$ if there is a $t\in U$ such that $tvx = tuy$. From the definition, it is clear that this relation is reflexive and symmetric. In order to see that this is also a transitive relation, assume that $(x,u) \sim (y,v)$ and $(y,v) \sim (z,w)$. So there exist $t,t^{\prime} \in U$ such that $tvx=tuy$ and $t^{\prime} wy = t^{\prime} vz$. Use $y$ as an intermediate to reach the equality $(t^{\prime}tv)wx = (t^{\prime}tv) uz$, which obviously implies that $(x,u) \sim (z,w)$. This means that $\sim$ is an equivalence relation.

Set $s/u$ for the equivalence class of $(s,u)$ under $\sim$ and let $S_U = \{s/u : s\in S, u\in U\}$. The operations ``$+$" and ``$\cdot$" on $S_U$ are defined as usual: $x/u + y/v = (xv + yu)/ uv$ and $x/u \cdot y/v = xy/uv$. It is a routine exercise to see that $(S_U,+,\cdot)$ is a semiring and $\gamma : S \rightarrow S_U$ defined by $\gamma(a) = a/1$ is a semiring homomorphism. We denote the extension of the ideal $I$ of $S$ in $S_U$ by $I\cdot S_U$ or simply $I S_U$.

Now let $I$ be an ideal of $S$ and define $I_U := \{x/u : x\in I, u\in U\}$. One can easily check that $I_U$ is an ideal of $S_U$. The set $I_U$ is called the localization of the ideal $I$ at $U$. We collect the basic properties of localization of ideals in the following:

\begin{proposition}
	
	Let $S$ be a semiring, $U$ an MC-set and $I,I^{\prime}$ ideals of $S$. Then the following statements hold:
	
	\begin{enumerate}
		
		\item $I_U = I\cdot S_U$
		\item If $I \subseteq I^{\prime}$ then $I_U \subseteq I^{\prime}_U$.
		\item $(I+I^{\prime})_U = I_U + I^{\prime}_U$.
		\item $(I \cap I^{\prime})_U = I_U \cap I^{\prime}_U$.
		
	\end{enumerate}
	
	\begin{proof}
		The proof of the statements (2), (3) and (4) is straightforward. We only prove (1) as an example: Let $S$ be a semiring, $U$ an MC-set and $I$ an ideal of $S$ and define $\gamma : S \rightarrow S_U$ by $\gamma(x) = x/1$. Obviously, if we take $x\in I$ and $u\in U$, then $x/u = (x/1) (1/u) = \gamma(x) \cdot (1/u) = x\cdot (1/u)$, which means that $x/u \in I\cdot S_U$. On the other hand, any element of $I\cdot S_U$ is of the form $\sum^n_{i=1} (x_i)(s_i / u_i)$. Consider the following calculation:
		
		$\sum^n_{i=1} (x_i)\cdot (s_i / u_i) = \sum^n_{i=1} (x_i / 1)(s_i / u_i) =\sum^n_{i=1}  (x_i s_i)/u_i = (\sum^n_{i=1} s_i a_i x_i)/u$\\ where $u = \prod u_i$ and $a_i = u_1 u_2 \cdots \widehat{u_i} \cdots u_n$.
		
		This shows that $\sum^n_{i=1} (x_i)(s_i / u_i) \in I_U$ and that is the proof of what it was claimed.
	\end{proof}
	
\end{proposition}

Let, for the moment, $J$ be an ideal of $S_U$ and define $\gamma : S \rightarrow S_U$ by $\gamma(a) = a/1$. We know that the contraction of the ideal $J$, i.e., $I_J = \gamma^{-1}(J) = \{x : x\in S , x/1 \in J\}$ is an ideal of $S$.

\begin{proposition}
	
	\label{extendedideal}
	
	Let $S$ be a semiring and $U$ an MC-set. Then every ideal of $S_U$ is an extended ideal.
	
	\begin{proof}
		Let $J$ be an ideal of $S_U$ and let $s/u \in J$. It is clear that this implies $s/1 \in J$ and therefore $s \in J^c$, which implies that $s/u \in J^{ce}$. But in general we know that $J \supseteq J^{ce}$. So $J = J^{ce}$, which means that every ideal of $S_U$ is an extended ideal.
	\end{proof}
	
\end{proposition}

\begin{remark}
	Let $S$ be a semiring and $I$ an ideal of $S$. The equality $I^{ec} = I$ is not always true even in commutative ring theory (Cf. \cite[Remark 5.27]{Sharp2000}). By the way, prime ideals behave much better as Theorem \ref{primesoflocalization} will show us:
	
\end{remark}

\begin{theorem}
	
	\label{primesoflocalization}
	
	Let $S$ be a semiring and $U$ an MC-set in $S$. Then the prime ideals of $S_U$ are in one-to-one correspondence with the prime ideals of $S$ disjoint from $U$.
	
	\begin{proof}
		Take $P$ to be a prime ideal of $S$, disjoint from $U$ and define $P^e= P_U = \{p/u : p\in P, u\in U\}$. It is easy to check that $P_U$ is a proper ideal of $S_U$. Now we prove that $P_U$ is a prime ideal of $S_U$. Take $x/u , y/v \in S_U$ such that $(x/u)(y/v) \in P_U$. This means that there are some $p\in P$ and $w\in U$ such that $xy/uv = p/w$ and so there is a $t\in U$ such that $twxy = tuvp \in P$. But $tw\in U$ and $U$ is disjoint from $P$ so $tw\notin P$ and therefore $xy\in P$ which implies that either $x\in P$ or $y\in P$. This shows us that either $x/u \in P_U$ or $y/v \in P_U$.
		
		It is easy to see that the map $P \mapsto P^e$ on the set of all prime ideals disjoint from $U$ is one-to-one. In order to prove that this map is onto, we must take a prime ideal $Q$ of $S_U$ and find a prime ideal $P_1$ disjoint from $U$ such that $P^e_1 = Q$. We set $P_1 = Q^c = \{x : x\in S, x/1\in Q\}$. It is easy to check that $P_1$ is a (proper) prime ideal $S$ disjoint from $U$. At last $P^e_1 = Q^{ce} = Q$, since in general we know that $J^{ce} = J$ and this finishes the proof.
	\end{proof}
	
\end{theorem}

Examples of MC-sets include the set of multiplicatively cancelable elements $\MC(S)$ of $S$ and $W = S-P$, where $P$ is a prime ideal of $S$. The case $W=S-P$ is of special interest and the reason is that the set $PS_W = \{x/w : x\in P, w\in S-P\}$ is the only maximal ideal of the semiring $S_W$ and the proof of our claim is as follows:

If $x/w \in S_W - PS_W$, then $x\notin P$, which means that $x/w$ is an invertible element of $S_W$. On the other hand, if $x/w$ is an invertible element of $S_W$, then $x\notin P$, which means that $x/w \in S_W - PS_W$. Therefore by Corollary \ref{localsemiring1}, $(S_W, PS_W)$ is a local semiring. The local semiring $S_W$ is usually denoted by $S_P$ and its unique maximal ideal $PS_W$ by $PS_P$ and the process of constructing $S_P$ from $S$ is called \emph{localization of $S$ at $P$}.

\begin{corollary}
	
	\label{primesoflocalizationatP}
	
	Let $S$ be a semiring and $P$ a prime ideal of $S$. Then the prime ideals of the local semiring $S_P$ are in one-to-one correspondence with the prime ideals of $S$ contained in $P$.
	
\end{corollary}

Let $S$ be a semiring and $(M,+,0)$ be a commutative monoid. The monoid $M$ is said to be an $S$-\emph{semimodule} if there is a function, called \emph{scalar product}, $\lambda: S \times M \rightarrow M$, defined by $\lambda (s,m)= s\cdot m$ such that the following conditions are satisfied:

\begin{enumerate}
	\item $s\cdot (m+n) = s\cdot m+s\cdot n$ for all $s\in S$ and $m,n \in M$;
	\item $(s+t)\cdot m = s\cdot m+t\cdot m$ and $(st)\cdot m = s\cdot (t\cdot m)$ for all $s,t\in S$ and $m\in M$;
	\item $s\cdot 0_M=0_M$ for all $s\in S$ and $0_S \cdot m = 0_M$ and $1_S\cdot m=m$ for all $m\in M$.
\end{enumerate}

A nonempty subset $N$ of an $S$-semimodule $M$ is said to be an $S$-subsemimodule of $M$ if $N$ is an $S$-semimodule itself.

Let $M$ be an $S$-semimodule. Then similar to semirings of fractions, one can see that the relation $\sim^{\prime}$ on $M \times U$, defined by $(m,u) \sim^{\prime} (n,v)$ if $tvm = tun$ for some $t\in U$ is an equivalence relation and if we put $m/u$ for the equivalence class of $(m,u)$ under $\sim^{\prime}$ and let $M_U = \{m/u : m\in M, u\in U\}$ and define addition ``$+$" and scalar product ``$\cdot$" as usual: $m/u + n/v = (vm + un)/ uv$ and $a/u \cdot m/v = am/uv$, then $M_U$ is an $S_U$-semimodule. Note that it is also possible to consider $M_U$ as an $S$-semimodule with the scalar product $s\cdot m/u = sm/u$ and therefore $\gamma : M \rightarrow M_U$ defined by $\gamma(m) = m/1$ is an $S$-semimodule homomorphism with this property that if $\gamma(m)=0$ then there exists a $t\in U$ such that $tm =0$.

\begin{proposition}
	
	Let $S$ be a semiring, $I$ an ideal of $S$, and $U$ an MC-set. Let $M$ be an $S$-semimodule and $K,L$ be $S$-subsemimodules of $M$. Then the following statements hold:
	
	\begin{enumerate}
		
		\item If $K \subseteq L$ then $K_U \subseteq L_U$.
		\item $(K+L)_U = K_U + L_U$.
		\item $(K \cap L)_U = K_U \cap L_U$.
		
		\item $(IL)_U = I_U L_U$.
		
	\end{enumerate}
	
	\begin{proof}
		Straightforward.
	\end{proof}
	
\end{proposition}

Let us recall that if $x$ is an element of an $S$-semimodule $M$, then the set $\Ann(x):=\{s\in S : s\cdot x =0\}$ is an ideal of $S$.

\begin{theorem}

	Let $M$ be an $S$-semimodule. Then the following statements are equivalent:
	
	\begin{enumerate}
		\item $M=0$
		\item $M_\mathfrak{p}=0$ for all $\mathfrak{p} \in \Spec(S)$
		\item $M_\mathfrak{m}=0$ for all $\mathfrak{m} \in \Max(S)$
	\end{enumerate}
	
	\begin{proof}
		It is clear that $(1) \Rightarrow (2) \Rightarrow (3)$. The proof of $(3) \Rightarrow (1)$ is as follows:
		
		Let $x\in M$. Consider the ideal $\Ann(x)$ of $S$. If $\Ann(x)=S$, then $x=0$. If $\Ann(x) \not= S$, then there is a maximal ideal $\mathfrak{m}$ of $S$ such that $\Ann(x) \subseteq \mathfrak{m}$. Since $x/1 = 0$ in $M_\mathfrak{m}$, there is an $s\in S-\mathfrak{m}$ such that $sx=0$, which means that $s\in \Ann(x)$, a contradiction. Therefore $M=0$ and the proof is complete.
	\end{proof}
	
\end{theorem}

\section{Primary ideals}\label{sec:primary}

Primary decomposition of ideals is an essential topic in traditional ideal theory in commutative rings. The main scope of this section is to investigate primary ideals of semirings. We also encourage the reader to see the paper by Lescot on prime and primary ideals of semirings \cite{Lescot2015}. 

Primary ideals for rings were introduced in commutative algebra by the German mathematician Emanuel Lasker (1868--1941) who was a student of David Hilbert (1862--1943). He was also a chess player and philosopher. We begin this section by defining primary ideals for semiring and then we bring their basic properties.

Let us recall that an ideal $Q$ of a semiring is called a primary ideal if $Q$ is a proper ideal of $S$ and $xy\in Q$ implies either $x\in Q$ or $y^n \in Q$ for some $n\in \mathbb N$ \cite[p. 92]{Golan1999(b)}.

\begin{proposition}
	Let $Q$ be a primary ideal of a semiring $S$. Then $\sqrt{Q}$ is the smallest prime ideal containing $Q$.
	
	\begin{proof}
		By Theorem \ref{krullintersectionprimes}, we only need to prove that $\sqrt{Q}$ is a prime ideal of $S$. Take $xy\in \sqrt{Q}$, then by definition of the radical of an ideal, there is an $m\in \mathbb N$ such that $x^m y^m \in Q$. Now by definition of primary ideals, either $x^m \in Q$ or $(x^m)^n \in Q$ for some $n\in \mathbb N$. This implies that either $x\in \sqrt{Q}$ or $y\in \sqrt{Q}$ and the proof is complete.
	\end{proof}
	
\end{proposition}

\begin{remark}
	
	If $Q$ is a primary ideal of $S$ and $P = \sqrt{Q}$, then $Q$ is said to be $P$-primary.
	
\end{remark}

\begin{proposition}
	
	\label{primary1}
	
	If $Q$ is an ideal of a semiring $S$ such that $\sqrt{Q}\in \Max(S)$, then $Q$ is a primary ideal of $S$. In particular, any power of a maximal ideal is a primary ideal.
	
	\begin{proof}
		Let $Q$ be an ideal of a semiring $S$ and $\sqrt{Q} = \mathfrak{m}$ such that $\mathfrak{m} \in \Max(S)$. Take $xy\in Q$ such that $y\notin \sqrt{Q}$. Since $\sqrt{Q} = \mathfrak{m}$ is a maximal ideal of $S$, $\sqrt{Q} + (y) = S$. This implies that $\sqrt{Q} + \sqrt{(y)} = S$ and therefore $Q + (y) = S$, which means that there are $a\in Q$ and $b\in S$ such that $a+by = 1$. From this we get that $ax + bxy = x$. Since $a, xy\in Q$, we get that $x\in Q$ and this finishes the proof.
	\end{proof}
	
\end{proposition}

Note that each prime ideal is also a primary ideal. Now we introduce another method for making primary ideals. Also, see Proposition \ref{primary3}.

\begin{proposition}
	
	\label{primary2}
	
	If all $Q_i$ for $1\leq i \leq n$ are $P$-primary, then so is $Q = \bigcap^n_{i=1} Q_i$.
	
	\begin{proof}
		Take $ab\in Q$, while $a\notin Q$. Then for any $1\leq i \leq n$, there is an $n_i$ such that $b^{n_i}\in Q_i$. This means that for any $1\leq i \leq n$, $b\in \sqrt{Q_i}$. Note that $\sqrt{Q} = \sqrt{\bigcap^n_{i=1} Q_i} = \bigcap^n_{i=1} \sqrt{Q_i} = P$. Hence, $b\in Q$. Q.E.D.
	\end{proof}
	
\end{proposition}

Let us define a new notation: For any ideal $I$ of $S$ and any element $x\in S$, we define $[I:x] := \{s\in S : sx\in I\}$.

\begin{proposition}
	
	\label{primary3}
	
	Let $S$ be a semiring, $x$ an element of $S$ and $Q$ be a $P$-primary ideal. The following statements hold:
	
	\begin{enumerate}
		
		\item If $x\in Q$, then $[Q:x]=S$.
		\item If $x\notin Q$, then $[Q:x]$ is a $P$-primary and $\sqrt{[Q:x]}=P$.
		\item If $x\notin P$, then $[Q:x]=Q$.
		
	\end{enumerate}
	
	\begin{proof}
		$(2)$: It is obvious that $Q \subseteq [Q:x]$. Now take $y\in [Q:x]$. So $xy\in Q$, which obviously implies that $y\in P$. This means that $Q \subseteq [Q:x] \subseteq P$ and therefore by taking radical, we get $\sqrt{[Q:x]}=P$. Now let $yz\in [Q:x]$. This means that $xyz\in Q$. Therefore if $y\notin Q$, we have $xz\in Q$, which means that $z\in[Q:x]$.
	\end{proof}
	
\end{proposition}

\section{Decomposition of Ideals}\label{sec:decompositions}

An ideal $I$ of a semiring $S$ is called irreducible if for any ideals $J$ and $K$ of $S$, $I = J \cap K$ implies that $I = J$ or $I = K$ \cite[p. 92]{Golan1999(b)}.

\begin{proposition}
	
	\label{irreducibleideals}
	
	Let $s$ be a nonzero element of a semiring $S$ and $I$ an ideal of $S$ not containing $s$. Then there exists an irreducible ideal $J$ of $S$ containing $I$ and not containing $s$.
	
	\begin{proof}
		Let $J_\alpha$ be a chain of ideals containing $I$ and not containing $s$. It is easy to check that $\bigcup_\alpha J_\alpha$ is also an ideal containing $I$ and not containing $s$. Therefore by Zorn's Lemma, we can find an ideal $J$ that is a maximal element of the set of all ideals of $S$ containing $I$ and not containing $s$. Imagine $J = K \cap L$, where $K$ and $L$ properly contain $J$. This implies that $a\in K$ and $a\in L$. But this means that $a\in K \cap L = J$, a contradiction. Therefore $J$ is irreducible.
	\end{proof}
	
\end{proposition}

\begin{proposition}
	
	If $I$ is a proper ideal of a semiring $S$, then $I$ is the intersection of all irreducible ideals of $S$ containing it.
	
	\begin{proof}
		Let $I$ be a proper ideal of $S$. This means that $1\notin I$. So by Proposition \ref{irreducibleideals}, there is an irreducible ideal containing $I$. Let $J$ be the intersection of all irreducible ideals of $S$ containing $I$. It is vivid that $I \subseteq J$. Our claim is that $I = J$. Suppose not. Then there is an element $s\in J-I$ and by Proposition \ref{irreducibleideals}, there is an irreducible ideal $K$ containing $I$ but not the element $s$ that is clearly a contradiction. Thus $I=J$ and the proof is complete.
	\end{proof}
	
\end{proposition}

\begin{proposition}
	Let $S$ be a Noetherian semiring. Then every ideal of $S$ can be represented as an intersection of a finite number of irreducible ideals of $S$.
	
	\begin{proof}
		Let $\mathcal{I}$ be the set of all ideals of $S$ which are not a finite intersection of irreducible ideals of $S$. We claim that $\mathcal{I} = \emptyset $. On the contrary, assume that $\mathcal{I} \neq \emptyset $. Since $S$ is Noetherian, $\mathcal{I}$ has a maximal element $I$. Since $I\in \mathcal{I}$, it is not a finite intersection of irreducible ideals of $S$. Especially it is not irreducible, which means that there are ideals $J$ and $K$ properly containing $I$ with $I = J \cap K$. Since $I$ is a maximal element of $\mathcal{I}$, $J,K\notin \mathcal{I}$. Therefore $J$ and $K$ are a finite intersection of irreducible ideals of $S$. But, then, $I = J \cap K$ is a finite intersection of irreducible ideals of $S$, a contradiction.
	\end{proof}
	
\end{proposition}

\begin{theorem}
	
	\label{subirreprimary}
	
	Let $S$ be a Noetherian semiring and $I$ a subtractive ideal of $S$. If $I$ is irreducible, then it is primary.
	
	\begin{proof}
		Let $I$ be a non-primary ideal of $S$. This means that there are $s,t \in S$ such that $st\in I$ but $t\notin I$ and $s^n \notin I$ for all $n\in \mathbb N$. Since $st\in I$, $t\in [I : s]$. But $t\notin I$. So $I \subset [I : s]$. Now by Proposition \ref{ideal:ideal}, we have that $[I : s^n] \subseteq [[I : s^n] : s] \subseteq [I : s^{n+1}]$, which gives us the following ascending chain of ideals:
		
		$ $
		
		\centerline{$I \subset [I : s] \subseteq \cdots \subseteq [I : s^n] \subseteq [I : s^{n+1}] \subseteq \cdots $.}
		
		$ $
		
		Since $S$ is Noetherian, this chain must stop somewhere, which means that there is some $m\in \mathbb N$ such that $  [I : s^m] = [I : s^{m+i}]$ for any $i\geq 0$. Our claim is that $I = [I : s^m] \cap (I+(s^m))$. Obviously, $[I : s^m]$ and $I+(s^m)$ contain $I$. Now let $x\in [I : s^m] \cap (I+(s^m))$. Since $x\in I+(s^m)$, there are some $y\in I, z \in S$ such that $x= y + zs^m$. But $x\in [I : s^m]$, which means that $ys^m+zs^{2m} = xs^m \in I$. Since $I$ is a subtractive ideal of $S$, we have $zs^{2m} \in I$, which means that $z\in [I : s^{2m}]$. But $[I : s^{2m}] = [I : s^m]$, so $zs^m \in I$ and this finally causes $x\in I$. This means that $I$ is reducible, the thing it was required to have shown.
	\end{proof}
	
\end{theorem}

Now we prove the so-called primary decomposition of ideals in semirings:

\begin{corollary}
	
	\label{primarydecomposition}
	
	Let $S$ be a subtractive Noetherian semiring. Then every ideal of $S$ can be represented as an intersection of a finite number of primary ideals of $S$.
\end{corollary}

A primary decomposition of an ideal $I$ of a semiring $S$ is a presentation of $I$ as a finite intersection of primary ideals of $S$ like the following:\\

\centerline{(PD) $I = \bigcap^n_{i=1} Q_i$, where $Q_i$ is a primary ideal of $S$ for any $1 \leq i \leq n$.\\}

$ $

If in addition, the prime ideals $P_i = \sqrt{Q_i}$ are all distinct and $Q_i \nsupseteq \bigcap_{j\neq i} Q_j$ for any $1 \leq i \leq n$, then it is said that the primary decomposition (PD) is minimal. Using Proposition \ref{primary2}, it is clear that any primary decomposition can be reduced to its minimal form. The prime ideals $P_i$ ($1 \leq i \leq n$) in the minimal decomposition of the ideal $I$ are said to belong to $I$. The minimal elements of the set of all primes belonging to $I$ are said to be minimal prime ideals belonging to $I$.

\begin{proposition}
	Let $S$ be a semiring and the ideal $I$ of $S$ possess a primary decomposition. Then the following statements hold:
	
	\begin{enumerate}
		
		\item Any prime ideal $P\supseteq I$ contains a minimal prime ideal belonging to $I$.
		
		\item The minimal prime ideals belonging to $I$ are precisely the minimal elements in the set of all prime ideals containing $I$.
		
	\end{enumerate}
	
	\begin{proof}
		$(1)$: Let $P\supseteq I$ be a prime ideal of $S$ and $I = \bigcap^n_{i=1} Q_i$ its minimal primary decomposition. Then $P = \sqrt{P} \supseteq \sqrt{I} = \sqrt{\bigcap^n_{i=1} Q_i} = \bigcap^n_{i=1} \sqrt{Q_i} = \bigcap^n_{i=1} P_i$. By Corollary \ref{primecontainingideals}, $P \supseteq P_i$ for some $i$. The statement $(2)$ is a direct consequence of $(1)$ and the proof is complete.
	\end{proof}
	
\end{proposition}

\section*{Acknowledgments}  The author is supported in part by the Department of Engineering Science at the Golpayegan University of Technology and wishes to thank the department for supplying all necessary facilities in pursuing this research. It is also a pleasure to thank Professor Dara Moazzami, the President of the Golpayegan University of Technology, for his help and encouragements.


\begin{thebibliography}{9}
	
\bibitem{AhsanMordesonShabir2012} Ahsan, J., Mordeson,  J.N., Shabir, M.: {\em Fuzzy Semirings with Applications to Automata Theory}, Springer, Berlin, 2012.
	
\bibitem{AtiyahMacDonald1969} Atiyah, M.F., MacDonald, I.G.: {\em An Introduction to Commutative Algebra}, Addison-Wesley, Reading, 1969.

\bibitem{Bistarelli2004} Bistarelli, S.: {\em Semirings for Soft Constraint Solving and Programming}, Springer-Verlag, Berlin, 2004.

\bibitem{Bourne1951} Bourne, S.: {\em The Jacobson radical of a semiring}, Proc. Nat. Acad. Sci. {\bf 37} (1951), 163--170.

\bibitem{Eilenberg1974} Eilenberg, S.: {\em Automata, Languages, and Machines}, Vol. A., Academic Press, New York, 1974.


\bibitem{GhalandarzadehNasehpourRazavi2017} Ghalandarzadeh, S., Nasehpour, P., Razavi, R.: {\em Invertible ideals and Gaussian semirings}, Arch. Math. Brno, {\bf 53}(3) (2017), 179--192. 

\bibitem{Gilmer1972} Gilmer, R.: {\em Multiplicative Ideal Theory}, Marcel Dekker, New York, 1972.

\bibitem{Glazek2002} G{\l}azek, K.: {\em A Guide to the Literature on Semirings and Their Applications in Mathematics and Information Sciences}, Kluwer, Dordrecht, 2002.

\bibitem{Golan1999(a)} Golan, J.S.: {\em Power Algebras over Semirings, With Applications in Mathematics and Computer Science}, Kluwer, Dordrecht, 1999.

\bibitem{Golan1999(b)} Golan, J.S.: {\em Semirings and their Applications}, Kluwer, Dordrecht, 1999.

\bibitem{Golan2003} Golan, J.S.: {\em Semirings and Affine Equations over Them: Theory and Applications}, Kluwer, Dordrecht, 2003.

\bibitem{Golan2005} Golan, J.S.: {\em Some recent applications of semiring theory}, International Conference on Algebra in Memory of Kostia Beider at National Cheng Kung University, Tainan, 2005.

\bibitem{GondranMinoux2008} Gondran, M., Minoux, M.: {\em Graphs, Dioids and Semirings}, Springer, New York, 2008.

\bibitem{HebischWeinert1998} Hebisch U., Weinert, H.J.: {\em Semirings, Algebraic Theory and Applications in Computer Science}, World Scientific, Singapore, 1998.

\bibitem{Hodges1979} Hodges, W.: {\em Krull implies Zorn}, J. London Math. Soc. (2), 19 (1979), 285--287.

\bibitem{Huckaba1988} Huckaba, J.A.: {\em Commutative Rings with Zero Divisors}, Marcel Dekker, New York, 1988.

\bibitem{Iseki1956} Is\'{e}ki, K.: {\em Ideal theory of semiring}, Proc. Japan. Acad., {\bf 32} (1956), 554-559.

\bibitem{KuichSalomaa1986} Kuich, W., Salomaa, A.: {\em Semirings, Automata, Languages}, Springer-Verlag, Berlin, 1986.

\bibitem{LaGrassa1995} LaGrassa, S.: {\em Semirings: Ideals and Polynomials}, Ph.D. Thesis, University of Iowa, 1995.

\bibitem{Lescot2015} Lescot, P.: {\em Prime and primary ideals in semirings}, Osaka J. Math. {\bf 52} (2015), 721--736.

\bibitem{Monk1969} Monk, J.D.: {\em Introduction to Set Theory}, McGraw-Hill, New York, 1969.

\bibitem{Nasehpour2016} Nasehpour, P.: {\em On the content of polynomials over semirings and its applications}, J. Algebra Appl., {\bf 15}, No. 5 (2016), 1650088 (32 pages).

\bibitem{NasehpourP2018} Nasehpour, P.: {\em On zero-divisors of semimodules and semialgebras}, arXiv:1702.00810 (2018).

\bibitem{Nasehpour2017} Nasehpour, P.: {\em Pseudocomplementation and minimal prime ideals in semirings}, arXiv:1703.08923 (2017), to appear in Algebra Univers.

\bibitem{NasehpourPP2018} Nasehpour, P.: {\em Some remarks on semirings and their ideals}, 	arXiv:1804.00593 (2018).

\bibitem{Nasehpour2018} Nasehpour, P.: {\em Valuation semirings},  J. Algebra. Appl, {\bf 16}, No. 4 (2018) 1850073 (23 pages).

\bibitem{NasehpourParvardi2018} Nasehpour, P., Parvardi, A.H. {\em Finitely additive, modular, and probability functions on pre-Semirings}, Comm. Algebra, {\bf 46}(7) (2018), 2968--2989.

\bibitem{Pilz1983} Pilz, G.: {\em Near-rings. The Theory and its Applications}, Rev. ed., North-Holland Mathematics Studies, {\bf 23}. Amsterdam - New York - Oxford: North-Holland Publishing Company, 1983. 

\bibitem{Scott1954} Scott, D.S.: {\em Prime ideal theorems for rings, lattices and Boolean algebras}, Bull. Amer. Math. Soc, {\bf 60} (1954), 390.

\bibitem{Sharp2000} Sharp, R.Y.: {\em Steps in Commutative Algebra}, 2d. edn., Cambridge Univ. Press, Cambridge, 2000.

\bibitem{Vandiver1934} Vandiver, H.S.: {\em Note on a simple type of algebra in which the cancellation law of addition does not hold}, Bull. Am. Math. Soc. {\bf40} (1934), 914--920. 

\bibitem{ZariskiSamuel1975} Zariski, O., Samuel, P.: {\em Commutative Algebra}, Vol. I. With the cooperation of I.S. Cohen. 2nd ed. (English), Graduate Texts in Mathematics. 28. New York - Heidelberg - Berlin: Springer-Verlag, 1975.

\bibitem{ZariskiSamuel1976} Zariski, O., Samuel, P.: {\em Commutative Algebra}, Vol. II. Reprint of the 1958-1960 Van Nostrand edition. (English), Graduate Texts in Mathematics. 29. New York - Heidelberg -Berlin: Springer-Verlag, 1976. 

\end{thebibliography}
\end{document}